\documentclass[12pt,A4]{article}
\usepackage{amssymb,amsmath, mathrsfs}

\newcommand{\pr}{{\bf P}}                     
\newcommand{\exn}{{\bf E\,}}                    
\newcommand{\ind}{{\bf 1}}                     
\newcommand{\ep}{\varepsilon}            
\newcommand{\f}{\varphi}            
\newcommand{\al}{\alpha}                 
\newcommand{\R}{\mathbb{R}}             
\newcommand{\Su}{S}             

\newtheorem{theo}{Theorem}
\newtheorem{rema}{Remark}

\newcommand{\zhe}[6]{#1 #2. {\em #3}, {\bf #5} (#4), #6.}
\newcommand{\kn}[4]{#1 {\em #2} (#3, #4).}
\newcommand{\inkn}[7]{#1 #2, in {\em #3} #4 (#5, #6),  #7.}

\title{Large deviation probabilities for sums of censored random variables with regularly varying distribution tails\footnote{Research supported by the University of Melbourne Faculty of Science Research Grant Support Scheme.}}
\date{}
\author{Aaron Chong\footnote{School of Mathematics and Statistics, The University of Melbourne, Parkville 3010, Australia;  e-mail: awychong@hotmail.com.}{\ } and Konstantin Borovkov\footnote{School of Mathematics and Statistics, The University of Melbourne, Parkville 3010, Australia; e-mail: borovkov@unimelb.edu.au (the corresponding author).}}

\begin{document}

\maketitle

\begin{abstract}
Let $ \xi_1, \xi_2,\ldots$ be a sequence of independent and identically distributed random variables with zero mean, finite second moment and   regularly varying right distribution tail. Motivated by a stop-loss insurance model, we consider  a threshold sequence $M_n\gg (n\ln n)^{1/2},$ $n\to \infty,$ and  establish the asymptotics of the probabilities of the  large deviations  of the form $ \sum_{j=1}^n(\xi_j \wedge M_n)>x$ in the whole spectrum of $x$-values in the region $O(M_n).$ The asymptotic representations for these probabilities obey the ``multiple large jumps principle" and have different forms in vicinities of the multiples $kM_n$ of the censoring threshold values, on the one hand, and inside intervals of the form $((k-1)M_n, kM_n),$ on the other. We show that there is a  ``smooth transition"
of these   representations   from one to the other
when the deviation~$x$ increases to a multiple of~$M_n$,  ``crosses" it and then moves away from it.



\end{abstract}

\section{Introduction and main results}

Let $\xi, \xi_1, \xi_2,\ldots$ be a sequence of independent identically distributed (i.i.d.) random variables, and $\Su_0:=0, $ $\Su_n:= \Su_{n-1}+ \xi_n$, $n\ge 1,$ be the random walk (RW) generated by this sequence. Such RWs are classical objects of probability theory and the study thereof is of great theoretical interest. At the same time, they are key components of numerous stochastic models used in various application areas such as statistics, risk theory and queueing, where one is usually dealing with situations where the number~$n$ of steps in the RW tends to be large.

The most fundamental result concerning the limiting behavior of the distribution of $S_n$ as $n\to\infty$ under   the classical condition  that
\begin{align}
\exn \xi =0, \quad \exn \xi^2 =1,
\label{01}
\end{align}
which is assumed to be met in this paper, is the central limit theorem stating that
\begin{align}
\sup_{x}|\pr (S_n >x)- \overline{\Phi}(xn^{-1/2}) |\to 0 \quad\mbox{as} \ n\to \infty,
\label{CLT}
\end{align}
where $\overline{\Phi} (z):= 1- \Phi(z)$ is the ``tail" of the standard normal distribution function~$\Phi$. However, for large deviations (LDs) $x =x_n\gg n^{1/2},$ the above result only states that $\pr (S_n >x)\to 0$ which is not good enough in many applications where one needs approximations with small relative errors.

It clearly follows from~\eqref{CLT}  that there is a sequence $\sigma_n\gg n^{1/2}$ such that
\begin{align}
\label{CLTr}
\pr (S_n >x)\sim \overline{\Phi} (xn^{-1/2}) \quad \mbox{as $n\to\infty$ uniformly in~$x\le  \sigma_n$},
\end{align}
where $\sim$ denotes asymptotic equivalence of functions: we write  $f\sim g$ iff $f/g\to 1$ in the respective limit. The threshold value $\sigma_n$ until which the normal approximation for the tail probability has vanishing relative error depends on the distribution of~$\xi$. Thus, under the additional moment  condition $\exn |\xi|^3<\infty$ ensuring that the Berry--Esseen convergence rate bound in the central limit theorem holds true, 
one can take $\sigma_n :=(1-\ep)(n \ln n)^{1/2}$ for any fixed $\ep>0$, whereas if $\f (\lambda):= \exn e^{\lambda \xi}<\infty,$   $|\lambda|< \lambda_0, $ for some $\lambda_0<\infty$ then one can take $\sigma_n := \ep_n n^{2/3}$ for any sequence $\ep_n\downarrow 0$, see e.g.\ Ch.~VIII in~\cite{Pe75} for more detail and further results.

The history of systematic work on specifically evaluating the (small) probabilities of LDs goes back to~\cite{Cr38}. The research  was first focussed on the case where the Cram\'er condition is met:
\begin{align}
\label{Cram}
\f (\lambda) <\infty\quad\mbox{for some }\ \lambda>0.
\end{align}
 The asymptotic behavior of  $\pr (S_n >x)$ in this case is formed, roughly speaking, by contributions made by all the jumps in the RW up to time~$n$. It is described by laws that are established mostly via analytical calculations and are determined by the properties of the moment generating function~$\f $ of~$\xi$. For more detail, see e.g.\ Ch.~VIII in~\cite{Pe75} and Ch.~9 in~\cite{Bo13}.

The situation changes dramatically when the right distribution tail of~$\xi$ is ``heavy", i.e., when the Cram\'er  condition~\eqref{Cram} is not met. It turns out that, under appropriate regularity conditions on the distribution tail $V(t):=\pr (\xi >t)$ as $t\to\infty$, the asymptotics of~$\pr (S_n >x)$ are governed by the ``single large jump principle" when~$x$ is ``large enough". This means that there is a single large jump in the trajectory $\{S_j\}_{0\le j\le n}$ up to time~$n$ that ``drives" the RW path above the high level~$x,$ after which the RW stays at roughly the same level until the time~$n,$ ensuring that the terminal value~$S_n$ also exceeds~$x$. In particular, under our assumptions~\eqref{01} and the additional   condition that the right tail of the jump distribution is a regularly varying function: 
\begin{align}
\label{regular_V}
V(t) = t^{-\al}L(t),\quad\mbox{where $\al>2,$ $L(t)$ is slowly varying at infinity}
\end{align}
(i.e., for any $y>0$, $L(ty)\sim L(t) $ as $t\to\infty$), the following result holds true. Set 
$$
s_n:= {((\al -2)n \ln n)}^{1/2},\quad n\ge  1.
$$ 
Then, for any fixed $c>1,$  as $n\to\infty$,
\begin{align}
\label{single_j}
\pr (S_n >x)\sim  n V(x) \quad \mbox{uniformly in  $x\ge c s_n$.}
\end{align}
It is easily seen that also  $\pr (\max_{j\le n}\xi_j >x)\sim nV(x),$ whereas $\pr \big(\sum_{j\le n}\ind_{\{\xi_j>x\}}\ge 2\big)\le \frac12 (n V(x))^2\ll nV(x)$ for $x\ge s_n$ (cf.~\eqref{nu>k}). This means that the probability in~\eqref{single_j} is essentially that of having a single jump of size greater than~$x$ in the first~$n$ steps of the RW.  Unlike the case of the light-tailed~$\xi$, the methods used to establish this kind of asymptotics are mostly probabilistic.

The first results in this direction were obtained in~\cite{Na65,He67,Na69,Na69a}; see monograph~\cite{BoBo08} for more references, bibliographic comments  and a systematic treatment of the relevant LD theory and results for regularly-varying and other classes of heavy-tailed distributions.

Note that relation~\eqref{single_j} for (any) {\em fixed~$n$ as $x\to\infty$} is essentially the subexponentiality property of the distribution of~$\xi$ (more precisely,  to be subexponential, the  distribution of a random variable that can take values of both signs must, in addition to  satisfying~\eqref{single_j}, be {\em long-tailed}: for any fixed $y\in \R,$ $V(t+y)\sim V(t)$ as $t\to\infty$, see e.g.\ Section~3.2 in~\cite{FoKoZa13}). Relation~\eqref{single_j}, however, holds when both~$x$  and~$n$ tend to infinity, and subexponentiality alone does not imply it in such a case. The ranges of the $n=n(x)$-values for which~\eqref{single_j} holds as~$x\to\infty$ for different classes of distributions were found   in Section~5.9 in \cite{BoBo08}; see also~\cite{BoBo13} for further results in that direction.

There is also a uniform representation combining both asymptotics~\eqref{CLTr}  and~\eqref{single_j} on the positive half-line  under assumptions~\eqref{01}, \eqref{regular_V} and the additional condition
\begin{align}
\label{condRoz}
\exn (\xi^2; |\xi|>t) = o (1/\ln t) \quad\mbox{as \ } t\to\infty.
\end{align}
Namely, setting
\[
 H_n(z):= \overline{\Phi} (zn^{-1/2})  +  n V(z)\ind_{\{z>n^{1/2}\}},\quad  z\in \R,
 \]
one has
\begin{align}
\label{uniRoz}
 \pr (S_n >x)\sim H_n(x) \quad\mbox{as \ }  n\to \infty \ \mbox{uniformly in~$x\in \R$}
\end{align}
(Corollary~7 in~\cite{Ro89}; it is also shown there that condition~\eqref{condRoz} is in fact  necessary for~\eqref{uniRoz} under~\eqref{01}, \eqref{regular_V}).
From the well-known asymptotics  $\overline{\Phi} (t)\sim   t^{-1} \phi (t)$ as $t\to\infty$ (Mills' ratio), where $\phi (t):=(2\pi )^{-1/2}  e^{-t^2/2}$ is the standard normal density, one can easily see that, in   approximation~\eqref{uniRoz}, the ``switch" in $H_n(x)$  from the normal term to the single jump asymptotics occurs just after the point~$x=s_n.$

Single jump asymptotics results akin to \eqref{single_j} were also obtained in the triangular array scheme, for sums $S_n:= \xi_{n,1}+\cdots+ \xi_{n,n}$, where $\xi_{n,1},\ldots, \xi_{n,n}$ are independent but not necessarily identically distributed random variables whose distributions can depend on~$n$ and should satisfy a suitable uniform regular variation condition  in place of~\eqref{regular_V}, see Ch.~13 in~\cite{BoBo08}.

In the present paper, we also deal with a triangular array scheme, but of a different kind. Namely, starting with our original sequence $\{\xi_j\}_{j\ge 1}$ satisfying conditions~\eqref{01} and~\eqref{regular_V}, we introduce a threshold sequence
\begin{align}
\label{M}
M_n \gg s_n \quad\mbox{as \ }  n\to \infty
\end{align}
and consider  the sums of ``censored" (or ``clipped") random variables $X_{n,j}:= \xi_j \wedge M_n,$ $j\ge 1$:
\[
Y_n:= \sum_{j=1}^n X_{n,j}, \quad n\ge 1,\quad  Y_0:=0.
\]
In the literature, the variables $X_{n,j}$ are sometimes referred to as ``truncated", but the latter term also applies to $ \xi_j \ind_{\{\xi_j\le  M_n\}},$ so we do not use it to avoid confusion. We study in this paper  the asymptotic behavior of the probabilities $\pr (Y_n >x)$ as $n\to\infty$ for different ranges of LDs~$x$.

The main motivation for considering the above model comes from actuarial practice. For instance, an insurance company can purchase an aggregate stop-loss insurance policy from a reinsurer. This policy protects the company against the combined cost of its customers' claims over a set amount within, say, a calendar month. Denote by~$\zeta_j$ the aggregate amount of all the claims made against the company during month~$j$ and assume that $\{\zeta_j\}_{j\ge 1}$ is a  sequence of i.i.d.\ random variables, with finite mean  $m:=\exn \zeta_1$. Suppose that the stop-loss insurance contract covering a time period of~$n\ge 1$ months specifies an amount~$N_n$ such that, each month, the reinsurer will reimburse the company  for the remainder of the that month's claims to be paid over   that  amount~$N_n$.
Assuming that our company receives constant monthly premium payments~$c$ and letting $\xi_j:=\zeta_j -m,$ $j\ge 1,$ $M_n:=N_n-m,$ we see that the claim surplus value for our company after~$n$ months is given by
\[
\sum_{j\le n}(\zeta_j\wedge N_n -c)= \sum_{j\le n}(\xi_j\wedge M_n  + m -c)=Y_n+n (m-c).
\]

Note that if the censoring threshold $M_n$ were {\em fixed} (i.e., independent of~$n$), then $\{X_{n,j}\}$ would be an i.i.d.\ sequence for which the Cram\'er condition is met (since $\exn e^{\lambda X_{n,j}}\le e^{\lambda M }<\infty$ for any $\lambda>0$), so that the Cram\'er LD theory would be applicable. Having $M_n\to\infty $ as~$n\to\infty$ makes the situation much more interesting and, at the same time, more relevant to practical problems as the stop-loss limits are usually rather high. Intuitively, one could expect that if~$M_n$ increases slowly with~$n$ then the LD probabilities behavior will be of the light-tail nature, whereas if~$M_n$ increases fast enough then the heavy-tailed asymptotics will be more relevant.

It turns out that, in this setup, instead of the single jump asymptotics~\eqref{single_j} valid for the original RW under assumptions~\eqref{01} and~\eqref{regular_V}, there emerge more complex multiple-jump asymptotics for $\pr (Y_n>x)$ depending on the range of~$x$ values. Our Theorems~\ref{thm_kM} and~\ref{thm_kinM} below provide complete description of their behavior, both when $x$ is ``passing" through a multiple $kM_n,$ $k=1,2,\ldots,$ of the censoring level and when $x$ is ``strictly inside" the interval between the consecutive multiples $(k-1)M_n$ and $kM_n$, with a continuous transition of one of these asymptotics to the other when the deviation~$x$ passes through  the vicinity of a multiple $kM_n$.

LD probabilities of the sums of censored i.i.d.\ random vectors $\boldsymbol{\xi},\boldsymbol{\xi}_1,\boldsymbol{\xi}_2,  \ldots \in \R^d$ with regularly varying tails (for the definition of multivariate regular variation, see e.g.~\cite{BaDaMi02,Re87}) were studied earlier  in~\cite{Ch12} for the following model discussed in~\cite{ChSa12}. In the latter paper, the authors considered the sums $\boldsymbol{Y}_{n }:=\sum_{j\le n}\boldsymbol{X}_{n,j}$ of isotropically  truncated   random vectors of the following form: for $M_n\to \infty,$
\begin{align}
\label{ChakrModel}
\boldsymbol{X}_{n,j} := \boldsymbol{\xi}_j \ind_{\{\|\boldsymbol{\xi}_j\|\le M_n \}}  +\frac{\boldsymbol{\xi}_j}{\|\boldsymbol{\xi}_j\|} (M_n+R_j)  \ind_{\|\{\boldsymbol{\xi}_j\| > M_n \} },
\end{align}
where $\{R_j\}_{j\ge 1}$ is a sequence of light-tailed non-negative i.i.d.\ random variables, independent of $\{\boldsymbol{\xi}_j\}_{j\ge 1}$. The authors called truncation ``soft" when $n\pr (\|\boldsymbol{\xi}_1\|>M_n) \to 0 $ as $n\to\infty,$  and ``hard" when $n\pr (\|\boldsymbol{\xi}_1\|>M_n) \to \infty.$

In the case when the tail exponent $\alpha$ of the $\boldsymbol{\xi}_j$'s is in the interval  $(0,2),$  it was shown in~\cite{ChSa12} that ``observations with softly truncated tails behave like heavy tailed random variables, while observations with hard truncated
tails behave like light tailed random variables" (meaning that the distributions of the appropriately scaled sums~$\boldsymbol{Y}_{n }$ would converge in the latter case to a normal law). The authors also considered statistical problems, on estimating the tail exponent $\alpha$ from truncated observations without knowing the truncation level~$M_n$, and on testing the hypothesis of the soft (correspondingly, hard) truncation regime against the appropriate alternative. That paper also contains references to some earlier work on RWs with censored jumps.

In~\cite{Ch12} the author proved, in the soft truncation regime (when $M_n\gg s_n$), the vague convergence of the measures  $\pr (M_n^{-1} \boldsymbol{Y}_n\in \cdot\,) (n\pr (\|\boldsymbol{\xi}\|>M_n)^{-k} $ on the spherical layer $\{\boldsymbol{x}\in \R^d: k-1<\|\boldsymbol{x}\|\le k\}$ as $n\to\infty $ (for $k>1$; the claim for $k=1$ is somewhat different). They also derived for that regime the decay rates in the ``boundary case" by finding, under suitable conditions, the weak limit (as $n\to\infty$) of the measures
\[
\pr (\|\boldsymbol{Y}_n\|>kM_n , \boldsymbol{Y}_n /\|\boldsymbol{Y}_n\|\in \cdot\,)(n\pr (\|\boldsymbol{\xi}\|>M_n))^{-k}
\]
on the unit sphere in~$\R^d$, $ k =1,2,\ldots $\      In the   hard truncation regime, it was shown in~\cite{Ch12} that the LD principle (referring to  the logarithmic asymptotics of the deviation probabilities) with ``speed" $n \pr (\|\boldsymbol{\xi}\|>M_n)$ holds for the scaled sequence $\{\boldsymbol{Y}_n/(nM_n\pr (\|\boldsymbol{\xi}\|>M_n))\}_{n\ge 1}$.

In~\cite{Ch17} these results were complemented by establishing in the univariate case the asymptotics of $\pr (Y'_n> k M_n),$ $k=1,2,\ldots,$   for a different truncation scheme, where one puts $Y'_n:= \sum_{j\le n }X'_{n,j}$ with $
X'_{n,j} := {\xi}_j \ind_{\{|{\xi}_j|\le M_n \}}$, and under the assumption that $M_n\gg n^{1/2 +\gamma}$ for some~$\gamma >0.$

The key difference of the model considered in the present paper from the   one in~\cite{Ch12,ChSa12} is that we assume that only the right tail of the jump distribution is regularly varying at infinity (whereas in the setting from~\cite{Ch12,ChSa12} one required ``two-sided regular variation" in the univariate case, including the possibility of the left tail vanishing faster that the regularly varying right one). In our setting, there are no ``extending" light-tailed random variables $R_j$ (cf.~\eqref{ChakrModel}), although one can include that term as well; under the assumption that it is bounded, one can show that all our results will still hold. The model in~\cite{Ch17}  imposes the regular variation condition on the right tail only, but it uses the ``true two-sided truncation" rather than censoring.

The key difference in the results obtained is that we establish the  asymptotics  for probabilities of the form $\pr (Y_n>x)$ for the whole spectrum of~$x$ values in the ``soft censoring  regime" rather than obtaining limiting laws for the scaled by~$M_n$ values of~$Y_n,$ as it was done in~\cite{Ch12,ChSa12}. The latter approach does not allow one to investigate what happens when~$x$ is in the vicinity of the multiples of~$M_n$. In this paper, we obtained fine asymptotic results for such~$x$'s  and   established  continuous transition of different kinds of approximations when~$x$ ``crosses" the boundary $kM_n$ between intervals of the form $((k-1)M_n, kM_n)$ and $(kM_n, (k+1)M_n),$ $k=1,2,\ldots$

Now we will state our main results. First of all, it turns out that, when the range of~$x$-values is ``noticeably" below~$M_n,$ there is no difference between the asymptotics of $\pr (Y_n>x)$ and those of $\pr(S_n>x).$  

\begin{theo}\label{thm_1M}
For any $c>1$ and sequence $\{d_n\}_{n\ge 1}$ such that $M_n - cs_n>d_n\gg n^{1/2}$ as $n\to \infty$, under assumptions~\eqref{01}, \eqref{regular_V} and~\eqref{M}, one has
\[
 \pr (Y_n >x)\sim n V(x)\quad\mbox{as \ }  n\to \infty
\]
uniformly in~$x\in (cs_n, M_n-d_n).$

Moreover, if we assume in addition that~\eqref{condRoz} is met then also
\[
 \pr (Y_n >x)\sim H_n(x) \quad\mbox{as \ }  n\to \infty
\]
uniformly in~$x\in (0, M_n-d_n).$
\end{theo}

Next we will turn to the case where $x$ is in the vicinity of a multiple of the censoring threshold~$M_n$. For convenience, we put
\[
\Pi_n :=nV(M_n), \quad n=1,2,\ldots 
\]
Since $M_n\gg s_n,$ one clearly has $\Pi_n\to 0$ as $n\to\infty.$

\begin{theo}\label{thm_kM}
Let conditions~\eqref{01}, \eqref{regular_V}, \eqref{condRoz} and~\eqref{M} be satisfied. Then, for any positive sequence $\ep_n \to 0$ as $n\to \infty$ and any fixed $k\ge 1,$ one has
\begin{align}
\pr (Y_n >x   ) \sim  \frac{ \Pi_n ^k}{k!} H_n (x-kM_n)  
\quad\mbox{as \ }  n\to \infty 
%
\label{0_sum}
\end{align}
uniformly in $x\in ((k-\ep_n)M_n, (k+\ep_n)M_n).$
\end{theo}

One can interpret~\eqref{0_sum} as follows. Roughly speaking, for $x$ close to $kM_n$, the most likely situation for the sum~$Y_n$ to exceed~$x$ is when 
\begin{itemize}
\item[(i)] exactly~$k$ of its terms are censored (each of them being then  equal to~$M_n$; note that $ \Pi_n^k/k!\sim {n\choose k} V^k (M_n)(1- V (M_n))^{n-k} $ as $ n\to \infty $), and
\item[(ii)]the sum of the remaining $n-k$ uncensored summands in $Y_n$ is greater than the difference $x-kM_n$, the probability of this event being approximately equal, by virtue of~\eqref{uniRoz}, to $H_{n-k} (x-kM_n)\sim H_n (x-kM_n) $ (see~\eqref{H=H}).
\end{itemize}    



To state our theorem concerning the asymptotics of $\pr (Y_n >x)$ when $x$ is ``well within" the interval between two consecutive multiples of~$M_n,$  we will need some further notations. For $k\ge 1 ,$    introduce simplices
\[
D_k (z) :=\Big\{\boldsymbol{t}= (t_1, \ldots, t_k)\in \R^k : \max_{1\le j\le k} t_j<1, t_1 + \cdots +t_k >z\Big\} , \quad  z\in (k-1,k),
\]
and  functions
\[
W_k (z):= \al^k\int_{
D_k (z)}   (t_1\cdots t_k)^{-\al -1} d\boldsymbol{t}, \quad  z\in (k-1,k).
\]
These functions can clearly be computed recursively:  setting $W_0(z)\equiv 1$, one has
\[
W_k (z) =\al\int_{z- (k-1)}^1 W_{k-1} (z-t) t^{-\al -1} dt, \quad  z\in (k-1,k), \ k\ge 1.
\]
In particular,
\begin{align}
  W_1 (z)  & = z^{-\al}-1,\quad z\in (0,1),
  \label{W1}
  \\
  W_2 (z) & = 1- (z-1)^{-\al} +\al z^{-2\al} {\rm B} (1-z^{-1}, z^{-1}; -\al, 1-\al),\quad z\in (1,2),
\notag
\end{align}
where ${\rm B}(z_1, z_2; a,b,)=\int_{z_1}^{z_2} t^{a-1} (1-t)^{b-1} dt,$ $0<z_1<z_2<1,$ $a,b\in \R,$  is the generalized incomplete beta function. For $k>2$, closed-form expressions for~$W_k$ become quite cumbersome.

\begin{theo}\label{thm_kinM}
Let conditions~\eqref{01}, \eqref{regular_V} and~\eqref{M} be satisfied. Then, for any fixed $k\ge 2$, there is a positive sequence $\{h_{k,n}\}_{n\ge 1}$ vanishing  slowly enough as $n\to \infty$ such that
\begin{align}
\pr (Y_n >x   ) \sim  \frac{  \Pi_n^k }{k!}  \sum_{j=0}^k {k\choose j} W_{k-j} (x/M_n -j), \quad n\to \infty,
\label{k_sum}
\end{align}
uniformly in $x\in ((k-1+h_{k,n})M_n, (k -h_{k,n})M_n).$
\end{theo}

\begin{rema}\label{rem_0}
{\rm Using the standard argument, one can see that, for any fixed $k_0\ge 2,$ the claims of both Theorems~\ref{thm_kM} and~\ref{thm_kinM} hold uniformly in $k\le k_0$ (with a common sequence $\{h_{k_0,n}\}$ in case of Theorem~\ref{thm_kinM}). Moreover, one can even assume here that $k_0\to\infty$ slowly enough.
}
\end{rema}

\begin{rema}\label{rem_1}
{\rm Since $W_0\equiv 1,$ the term with $j=k$ contributes $\Pi_n^k/k!$ to the right-hand side of~\eqref{k_sum}, whereas all the other terms in that expression are positive.
}
\end{rema}

\begin{figure}[ht]
\setlength{\unitlength}{1 mm}
\begin{center}
\begin{picture}(75, 63)(0,2)
%
%
\put(0,5){\vector(1,0){61}}
\put(5,0){\vector(0,1){60}}
\put(4,37){\line(1,0){52}}
\put(17,37){\line(0,1){20}}
\put(37,4){\line(0,1){52}}
\put(37,17){\line(1,0){20}}
{\small
\put(1,48){$x$}
\put(4,63){$\xi_2$}
\put(-1.5,36){$M_n$}
\put(2,1){$0$}
\put(35,0.5){$ M_n$}
\put(48,0.5){$x$}
\put(63,4){$\xi_1$}
\put(26,44){A}
\put(44,44){B}
\put(30,30){C}
\put(44,26){D}
} {\thicklines
\put(0,54){\line(1,-1){53}}
}
\end{picture}
\end{center}
\caption{\small An illustration to representation~\eqref{k_sum} in the case $k=2$ (so that $x\in (M_n, 2M_n)$), assuming that the two ``large" summands  are~$\xi_1$ and~$\xi_2$. The region labelled  with~B contributes to the term with~$j=2$ on the right hand side of~\eqref{k_sum} (both~$\xi_1$ and~$\xi_2$ are censored),   regions~A and~D contribute to the term with~$j=1$ ($\xi_2$ is censored in~A, $\xi_1$ is censored in~D), and   region~C corresponds to the term with~$j=0$ (none of the~$\xi_i$'s is censored). Roughly speaking, it is only when the ``large pair" $(\xi_1,\xi_2)$ is in one of the regions A--D that $Y_n$ will exceed~$x.$
}
\label{fig1}
\end{figure}
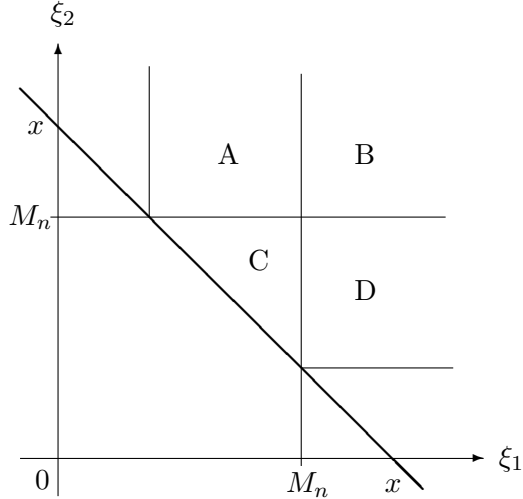

The meaning of the asymptotics in~\eqref{k_sum} can be loosely  explained as follows. For~$x$'s from the stated range, the sum~$Y_n$ exceeds~$x$ when it contains exactly~$k$ ``large  terms" (this refers, roughly  speaking, to the terms that are of the order of magnitude of~$M_n$). The $j$th term in the sum in~\eqref{k_sum} corresponds to the event that~$j$ of these~$k$ large terms were censored $\xi_i$'s (and so equal to~$M_n$ each) and the sum of the remaining $k-j$ large terms was big enough to ``bridge the remaining  gap" between~$x$ and~$jM_n$, the function~$W_{k-j}$ emerging from convolving $k-j$ power tail distributions approximating the conditional law of~$\xi/N$ given $ \xi>N\to\infty$ (cf.~\eqref{conv_cond}). In particular, the term with~$j=k$ in~\eqref{k_sum} corresponds to the event that all the~$k$ ``large summands" in~$Y_n$ were equal to~$M_n$, as in region~B in Fig.~\ref{fig1} illustrating the meaning of~\eqref{k_sum} in the case $k=2$, when there are four regions (A--D) contributing to the probability $\pr (Y_n>x)$ as explaind in the caption to the figure. In the case $k=3,$ one would have $2^3=8$ such regions, with three of them representing the event that only one of the three ``large" $\xi_i$'s is censored, a further three   corresponding to the event that two of these $\xi_i$'s are censored, with the two remaining regions having all the ``large" $\xi_i$'s censored in one and none in the other, and likewise for higher values of~$k.$

\begin{rema}\label{rem_3}
{\rm
Note that our theorems   provide together a complete coverage of the range of deviations that are $O(M_n)$ and establish a ``smooth transition" of the asymptotic representations from Theorems~\ref{thm_1M}--\ref{thm_kinM} from one to the other when the deviation~$x$ increases to a multiple of~$M_n$, ``crosses" it and then ``departs" from it.

Indeed, fix  any $k\ge 2$ (for simplicity of exposition; the same argument works for $k=1$ as well). Next set, say,  $\ep_n:=2 h_{k,n}$ in the statement of Theorem~\ref{thm_kM}  (assuming without loss of generality that $h_{k+1,n}= h_{k,n}$), so that there will be  overlaps of the form (see Fig.~\ref{fig2})
\[
I^{(n)}_-:= ((k -2h_{k,n})M_n,(k -h_{k,n})M_n), \quad I^{(n)}_+:= ((k +h_{k,n})M_n,(k + 2 h_{k,n})M_n)
\]
of the interval  where asymptotics~\eqref{0_sum} hold  in vicinity of~$kM_n,$ on the one hand, and the left- and right-adjacent intervals $ ((k-1+h_{k,n})M_n, (k -h_{k,n})M_n)$ and  $  ((k +h_{k,n})M_n, (k+1 -h_{k,n})M_n)$, respectively,  where representations of the form~\eqref{k_sum} hold, on the other  (in the case of~$I^{(n)}_+$, one must replace~$k$ with  $k+1$ in the expression on the right-hand side of~\eqref{k_sum}).

\begin{figure}[h]
\setlength{\unitlength}{1 mm}
\begin{center}
\begin{picture}(110, 10)(0,2)
{\thicklines
\put(0,3){\line(1,0){110}}
\multiput(15,2)(20,0){5}{\line(0,1){2}}
}
\put(15,5){$\overbrace{\hspace{20 mm}}^{\mbox{$I^{(n)}_-$}}$}
\put(75,5){$\overbrace{\hspace{20 mm}}^{\mbox{$I^{(n)}_+$}}$}
{\small \put(50,-2){$kM_n$}
\put(6,-2){  $(k-\ep_n)M_n$}
\put(86,-2){$(k+\ep_n)M_n$}
}
\end{picture}
\end{center}
\caption{\small The ``overlap intervals" $I^{(n)}_\pm$.}
\label{fig2}
\end{figure}
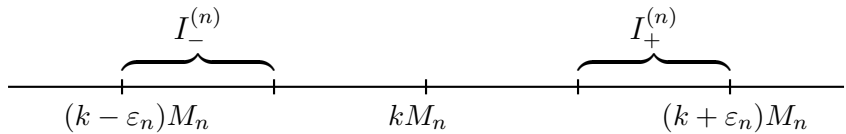

It will be seen in the proof of Theorem~\ref{thm_kinM} that $h_{k,n}$ should be chosen so that, in particular,
\begin{align}\label{h_prop}
 h_{k,n} M_n \gg s_n, \quad L(h_{k,n}M_n)\sim L(M_n)\quad \mbox{as $n\to\infty$}
\end{align}
(see~\eqref{hML} and the text prior to it).

In the ``left overlap" zone $I^{(n)}_-$,  all the terms in the sum in~\eqref{k_sum}   with $j<k$  are vanishing as $n\to\infty$  since $W_m(z)\downarrow 0$ as $z\uparrow m,$ $m\ge 1$ (as clearly $ D_m(z)\downarrow \varnothing$ then). As $W_0\equiv 1$, in the interval~$I^{(n)}_-$ the asymptotics turn  into just $\Pi_n^k/k!$. But this is exactly  what~\eqref{0_sum} becomes in~$I^{(n)}_-$  as well since, due to the first  relation in~\eqref{h_prop},
\[
1\ge \inf_{x\in I^{(n)}_-}H_n (x- kM_n)\ge   H_n (-h_{k,n}M_n)\to 1.
\]
In the ``right overlap" zone $I^{(n)}_+$, also due to the first relation in~\eqref{h_prop}, the right-hand side of~\eqref{0_sum} is equivalent to \begin{align}\label{right_zone}
nV(x-kM_n)\Pi_n^k/k!.
\end{align}
Now  we have to compare this with  the behavior in~$I^{(n)}_+$ of the expression on the right-hand side of~\eqref{k_sum} with~$k$ replaced by $k+1$. Consider the contributions of the terms with different values of $j=0,1,\ldots, k+1$ to the sum in that expression.

{\em Case $j=k+1$}.  As $W_0\equiv 1,$ the term with $j=k+1$ contributes $\Pi_n^{k+1}/(k+1)!$, which is negligibly small compared to~\eqref{right_zone} since $x-kM_n<2h_{k,n}M_n$ in~$I^{(n)}_+ $ and hence, using the second relation in~\eqref{h_prop}, one has
\[
\inf_{x\in I^{(n)}_+} V(x-kM_n) \ge V(2h_{k,n}M_n) = \frac{L(2h_{k,n}M_n)}{(2h_{k,n}M_n)^\alpha}
 \sim \frac{V( M_n)}{(2h_{k,n} )^\alpha}\gg V(M_n).
\]

{\em Case $j=k$}.  Due to~\eqref{W1} and the fact that $x/M_n -k\le 2 h_{k,n}\to 0$ in $I^{(n)}_+$, the contribution of the term with  $j=k$ to the right-hand side of~\eqref{k_sum} (with $k$ replaced by $k+1$) in that zone   is given by
\[
\frac{\Pi_n^{k+1}}{(k+1)!}{k+1 \choose k}W_1 (x/M_n -k)\sim
\frac{\Pi_n^{k+1}M_n^\alpha }{k!(x-kM_n)^\alpha}
 = \frac{\Pi_n^{k }}{k!}n (x-kM_n)^{-\alpha}L(M_n),
\]
which is asymptotically equivalent to~\eqref{right_zone} as  $  L(M_n)\sim L(x-kM_n)$ in~$I^{(n)}_+$ due to the second relation in~\eqref{h_prop}.

{\em Case $j<k$}. We will  show that,  for $x\in I^{(n)}_+$, the contribution  of each  the terms  with   $j<k$ to the right-hand side of~\eqref{k_sum} (with~$k$ replaced by $k +1$) is negligibly small compared to~\eqref{right_zone}.

To this end, for any~$m:=k+1-j\ge 2,$ introduce half-spaces
\[
E_{m,i}:=\{\boldsymbol{t}\in \R^m: t_i <1/3\},\quad
i\in [m]:=\{1, \ldots, m\},
\]
and note that, for any fixed $z\in (m-1,m),$  the sets $D_m (z)\cap  E_{m,i}$, $i\in [m],$ are   disjoint (indeed, if, say, $\boldsymbol{t} \in E_{m,1} E_{m,2} $ then  $\sum_{i=1}^m t_i <2/3+ (m-2)<m-1<z,$ so that $\boldsymbol{t}\not\in D_m (z)$). Therefore,  setting $E':= \bigcap_{i=1}^m E_{m,i}^c$,  due to symmetry one has
\begin{align*}
W_m(z) & = \alpha \int_{D_m (z) E' }  (t_1\cdots t_m)^{-\al -1} d\boldsymbol{t}
 + m\alpha \int_{D_m (z) E_{m,m} }  (t_1\cdots t_m)^{-\al -1} d\boldsymbol{t}.
\end{align*}
Denoting by $ {\rm vol}_j(\cdot ) $ the $j$-dimensional volume measure, we note that  the first integral clearly does not exceed
\[
 \sup_{\boldsymbol{t}\in D_m (z) E'} (t_1\cdots t_m)^{-\al -1} \cdot  {\rm vol}_m (D_m (z)) \le 3^{m(\al+1)} /m!
\]
and, as $ {\rm vol}_{m-1} \big(D_m (z) \{\boldsymbol{t}\in \R^m:t_m=t\}\big) < t^{m-1}/(m-1)!,$ $z\in (m-1,m),$  $t\in (0,1),$  the second integral is less than
\begin{align*}
 \sup_{\boldsymbol{t}\in D_m (z) E_{m,m}} (t_1\cdots t_{m-1})^{-\al -1}
 \int_{z- (m-1)}^{1/3} \frac{t^{m-\alpha -2} dt}{(m-1)!}
 =O((z-m+1)^{ -\alpha +1 })
  \end{align*}
as $z\downarrow m-1$ (note that $t^{m-\alpha -2}\le t^{-\alpha}$ since $m\ge 2$). Thus we showed that
\[
W_{k+1-j}(z) = O\big((z-k+j)^{-\alpha+1}\big), \quad z\downarrow k-j\ge 1.
\]

Therefore, for $x\in I^{(n)}_+,$ the contribution of the $j$th term in the sum in~\eqref{k_sum} (with~$k$ replaced by $k+1$) to the right-hand side of that relation is
\begin{align*}
\frac{\Pi_n^{k+1}}{(k+1)!}{k+1 \choose j}W_{k+1-j} (x/M_n-j)
 & =
\Pi_n^{k+1 } \cdot O(  (x/M_n-k)^{-\alpha+1}).
  \end{align*}
This is negligibly small compared to~\eqref{right_zone} since, due to the second relation in~\eqref{h_prop} (recall that   $ h_{k,n}  M_n <x-k M_n< 2h_{k,n} M_n$ for $x\in I^{(n)}_+$), one has
\begin{align*}
 \Pi_n (x/M_n-k)^{-\alpha+1}
 & = n L(M_n)   (x-kM_n)^{-\alpha}  (x/M_n-k)
 \\
 & \sim n V(x-kM_n)   (x/M_n-k)
 \\
 & \le 2h_{k,n}n V(x-kM_n) =o(n V(x-kM_n) ).
 \end{align*}
Thus, we showed that the asymptotics from our Theorems~\ref{thm_kM} and~\ref{thm_kinM} ``merge" in the right overlap zones~$I^{(n)}_+$ as well.
}
\end{rema}

\section{Proofs}

For a real sequence $\{a_j\}_{j\ge 1}$ and  $k\in [n],$  $n\ge 1,$ put
\[
\overline{a}_n:=\max_{j\in [n]}a_j, \quad \overline{a}_{n,k}:=\max_{n-k<j\le n}a_j, \quad
\underline{a}_n:=\min_{j\in [n]}a_j, \quad \underline{a}_{n,k}:=\min_{n-k<j\le n}a_j.
\]

{\em Proof of Theorem~\ref{thm_1M}.} For $y\in\R$, denote by
\[
\nu_n (y):=\sum_{j=1}^n \ind_{\{\xi_j> y\}}
\]
the number of jumps in the first~$n\ge 1$ steps of the RW $\{S_j \}_{j\ge 1} $ that exceed~$y.$
Note that, for any fixed $k=1,2,\ldots, $ one has, for $n\ge k,$
\begin{align}
\pr (  \nu_n (y)\ge  k  )
& \le \sum_{J\subseteq [n],\, |J|=k }\pr (\xi_i >y,   i\in J)
\notag
\\ &
={n \choose k }  V^{k } (y)\le \frac{(nV(y))^{k } }{k!}.
\label{nu>k}
\end{align}

As clearly $Y_n \le S_n $ always and $Y_n = S_n $ when $\nu_n (M_n)=0,$ one gets
\begin{align*}
  0 & \le \pr (S_n >x) - \pr (Y_n >x) = \pr (Y_n\le x< S_n)
  \\ & =   \pr (Y_n\le x< S_n, \nu_n (M_n)=1)  + \pr (Y_n\le x< S_n, \nu_n (M_n)\ge 2) .
\end{align*}
The last probability    does not exceed $\pr (\nu_n (M_n)\ge 2)=O(\Pi_n^2)$ by~\eqref{nu>k}, whereas the second last one, due to the $\xi_j$'s being i.i.d., equals
\begin{align*}
n\pr (S_{n-1} +M_n &\le x<S_n, \xi_n > M_n, \overline{\xi}_{n-1}\le M_n)
 \\
 & \le  n\pr (S_{n-1} +M_n \le x , \xi_n > M_n)
 \\
 & = \pr (S_{n-1}   \le x -M_n) \Pi_n \le (n-1) d_n^{-2}  \Pi_n =o(\Pi_n),
\end{align*}
where we used independence of $S_{n-1}$ and $\xi_n$, our assumption that $x <M-d_n,$  and Chebyshev's inequality to get the last two relations. Now the assertions of  Theorem~\ref{thm_1M} immediately follow from~\eqref{single_j} and~\eqref{uniRoz}, as $o(\Pi_n)$ is negligibly small compared to the main terms in these representations   in the respective ranges of $x$-values.

\hfill$\Box$

\medskip

{\em Proof of Theorem~\ref{thm_kM}.} Choosing $y:=M_n/(k+2)$, one has
\begin{align}
\label{1_sum}
\pr (Y_n >x )
& = \sum_{j=0}^{k-1} \pr (Y_n >x ,\nu_n (y)=j)
\notag \\
&  + \pr (Y_n >x , \nu_n (y)=k)  + \pr (Y_n >x , \nu_n (y)> k  ).
\end{align}
Here, for $0\le j<k ,$ 
one has
\begin{align}
 \pr (Y_n >x ,\nu_n (y)=j)
& = \sum_{J\subseteq [n], \, |J|=j} \pr (Y_n >x  , \xi_i>y, i\in J, \xi_m\le y,  m\in [n]\setminus J  )
\notag\\
& \le  {n\choose j}  \pr (S_{n-j} >x-jM_n, \overline{\xi}_{n-j} \le y , \underline{\xi}_{n,j} >y)
\notag
 \\
  & \le  \frac{n^j}{j!} V^j (y) \pr \big(S_{n-j} >(1-\ep_n)M_n, \overline{\xi}_{n-j} \le y  \big)
 \label{11_sum}
\end{align}
since $Y_n\le S_{n-j}+jM_n$, 
\[
x-jM_n \ge x-(k-1)M_n > (1-\ep_n)M_n\quad\mbox{for $x> (k-\ep_n)M_n,$}
\]
the $\xi_i$'s are independent and $\pr (\underline{\xi}_{n,j} >y)=\prod_{i=1}^j\pr (\xi_i>y)=V^j(y).$  Furthermore,  $nV(y)=o(1)$ as $y\gg s_n$ due to~\eqref{M}, so that $ n^j V^j (y)= O(1)$ for $0\le j<k.$ Finally, as $(1-\ep_n)M_n/y= (k+2)(1+o(1))$,  by Corollary~4.1.3 from~\cite{BoBo08}  one has, for any $\delta \in (0,1),$
\[
\pr \big(S _{n-j} >(1-\ep_n)M_n, \overline{\xi}_{n-j} \le y\big)=O\big(\Pi_n^{k+2 -\delta}\big) \quad   \mbox{as} \quad n\to\infty.
\]
We conclude that the left-hand side of~\eqref{11_sum} and hence the first term on the right-hand side of~\eqref{1_sum} is $o\big(\Pi_n^{k+1}\big).$

The last term on the right-hand side of~\eqref{1_sum} is $O\big(\Pi_n^{k+1}\big)$ by~\eqref{nu>k}.

Hence it remains to evaluate the middle term on the right-hand side of~\eqref{1_sum}. First let $\ep_n':=\ep_n + \eta_n,$ where we choose $\eta_n\downarrow 0$ so slowly as $n\to \infty$ that
\begin{align}
{\rm (i)}&~\eta_n M_n \gg s_n,
\notag
\\
 {\rm (ii)}&~V(\eta_n M_n)\sim \eta_n^{-\alpha} V(M_n),\label{eta}
 \\%
 {\rm (iii)}&~V(\ep_n M_n)\gg \eta_n^{-\alpha}  V(M_n) ,
\notag
\end{align}
which is always possible in view of $M_n\gg s_n$  (for~(i)) and the key properties of slowly  varying functions (see e.g.\ Theorem~1.1.2 and the remark after it in~\cite{BoBo08} for~(ii), and Theorem~1.1.4(iii) from~\cite{BoBo08} or Theorem~1.5.6 in~\cite{BiGoTe87} for~(iii) where it suffices to take $\eta_n\ge \ep_n^{1-\delta }$ for some $\delta>0$).

Set  
\[
A=A(n,k,x,y):= \{Y_n >x, \overline{\xi}_{n-k} \le y <\underline{\xi}_{n,k}\}.
\]
Due to the i.i.d.\ assumption on the $\xi_j$'s, one has
\begin{align}
\pr (Y_n >x , \nu_n (y)=k) ={n \choose k }\pr (A) ={n \choose k }(P_1 + P_2 +P_3),
\label{2b_sum}
\end{align}
where each of the three terms
\begin{align*}
P_1 &:= \pr \big(A;   \underline{\xi}_{n,k} \le (1-\ep'_n) M_n\big),
\\
P_2 &:=\pr \big(A; (1-\ep'_n) M_n <\underline{\xi}_{n,k} \le M_n\big),
\\
P_3 &:= \pr \big(A;    \underline{\xi}_{n,k} > M_n\big)
\end{align*}
will have to be evaluated in a different way.

First, as $\sum_{j=0}^{k-1} (\xi_{n-j} \wedge  M_n)\le (k-\ep'_n)M_n $ on the event in the definition of~$P_1$, and $x-kM_n> -\ep_n M_n,$ by independence of the $\xi_j$'s we get
\begin{align}
P_1 & \le \pr (S_{n-k}  + (k-\ep_n')M_n>x, \underline{\xi}_{n,k}>y )
\notag\\
&
 =\pr (S_{n-k}  >   (x-kM_n) +\ep'_nM_n) \pr^k ( \xi >y )
 \notag\\
 & \le \pr (S_{n-k}  >   \eta_n M_n) V^k(y)
  \notag\\
 & \sim n V( \eta_n M_n) V^k(y)= O\big(\eta_n^{-\al } n V^{k+1}(M_n)\big),
\label{2a_sum}
\end{align}
where the ``single jump  asymptotics" for $\pr (S_{n-k}  >   \eta_n M_n)$ hold true since $\eta_nM_n\gg s_{n-k}$ from~(i) in~\eqref{eta}, whereas the last relation in~\eqref{2a_sum} follows from~(ii) in~\eqref{eta}.  We conclude that, in view of~(iii) in~\eqref{eta}, the contribution of the term with~$P_1$ on the right-hand side of~\eqref{2b_sum} is
\begin{align}
{n \choose k } P_1= o\big(nV(\ep_n M_n)\cdot  \Pi^k\big).
\label{P1}
\end{align}

Second, again using symmetry and  the independence of the $\xi_j$'s,
\begin{align}
P_2 & \le k \pr \big(A; \xi_n\in ((1-\ep'_n)M_n , M_n]\big)
 \notag\\
 &\le k \pr \big( S_{n-k}+kM_n >x,    \underline{\xi}_{n-1,k-1}>(1-\ep'_n)M_n, \xi_n\in ((1-\ep'_n)M_n , M_n]\big)
 \notag\\
 & = k \pr  ( S_{n-k}>x-kM_n) V^{k-1}((1-\ep'_n)M_n) \big(V ((1-\ep'_n)M_n)-V ( M_n)\big)
 \notag\\
 & \sim k H_{n-k } (x-kM_n) V^{k-1}( M_n) \big(V ((1-\ep'_n)M_n)-V ( M_n)\big)
\label{2c_sum}
\end{align}
uniformly in $x\in ((k-\ep_n)M_n, (k+\ep_n)M_n),$ where the last relation in~\eqref{2c_sum} follows from 
the uniform representation~\eqref{uniRoz}. 

Next  we  note that
\begin{align}\label{H=H}
H_{n-k} (x-kM_n)\sim H_n (x-kM_n), \quad x\in ((k-\ep_n)M_n, (k+\ep_n)M_n).
\end{align}
To show this, it suffices to demonstrate that, as $n\to\infty,$
\begin{align}\label{Phii}
\overline{\Phi} ((x-kM_n)(n-k)^{-1/2})\sim \overline{\Phi} ((x-kM_n)n^{-1/2}),
\quad
x-kM_n\ll n,
\end{align}
as the equivalence of the ``single large jump'' terms in the representations for $H_{n-k}$ and $H_n$ in the range of the $x$-values where these terms dominate is obvious. The equivalence in~\eqref{Phii} clearly holds for $x-kM_n\le 0,$  whereas for   $0<x-k M_n\ll n $ this relation follows from the following bounds: setting 
$a:=(x-kM_n)n^{-1/2}<(x-kM_n)(n-k)^{-1/2}=:b,$ one has
\begin{align*}
\bigg| \frac{\overline{\Phi} (b)}{\overline{\Phi} (a)}-1\bigg|
= \frac{1}{\overline{\Phi} (a)}\int_a^b \phi (z) dz
\le \frac{(b-a)\phi (a)}{\overline{\Phi} (a)}
\le\bigg( \frac{b}{a}-1\bigg)(a^2+1)
\end{align*}
by the well-known inequality $\overline{\Phi} (z)/\phi(z) \ge z/(z^2+1),$ $z\ge 0,$   for Mills' ratio  from~\cite{Go41}. Since  $\frac{b}{a}-1=(\frac{n}{n-k})^{1/2}-1\sim \frac{ k}{2n} $  and $a^2=o(n),$  the desired equivalence follows.

As clearly $V ((1-\ep )M_n)-V ( M_n)=V(M_n) ((1-\ep)^{-\al} L((1-\ep)M_n)/L(M_n)-1)=o(V(M_n))$ as $\ep\to 0 $
from the uniform convergence theorem for slowly varying functions (see e.g.\ Theorem~1.1.2 in~\cite{BoBo08}  or Theorem~1.1.2 in~\cite{BiGoTe87}), we conclude from~\eqref{2c_sum} that
\begin{align}
{n \choose k } P_2= o\big(H_n (x-kM_n) \Pi_n^k\big).
\label{P2}
\end{align}

Third,
\begin{align}
P_3 & =  \pr \big(Y_n >x, \overline{\xi}_{n-k} \le y, \underline{\xi}_{n,k}>M_n\big)
 \notag\\
 & =\pr \big( S_{n-k}>x-kM_n , \overline{\xi}_{n-k} \le y\big)V^k(M_n)
 \notag\\
 & =\big[\pr \big( S_{n-k}>x-kM_n \big)
 - \pr \big( S_{n-k}>x-kM_n , \overline{\xi}_{n-k} > y\big)\big] V^k(M_n).
\label{2d_sum}
\end{align}
Here $\pr  ( S_{n-k}>x-kM_n  ) \sim H_n (x-kM_n)$ by~\eqref{uniRoz} and~\eqref{H=H}, and
\[
\pr  ( S_{n-k}>x-kM_n , \overline{\xi}_{n-k} > y )\le \pr  (   \overline{\xi}_{n-k} > y )\le (n-k)V(y)=O(\Pi_n),\]
implying that
\[
{n \choose k } P_3 = {n \choose k } H_n (x-kM_n)V^k(M_n) +O\big(\Pi_n^{k+1}\big).
\]
Since ${n \choose k }\sim n^k/k! $ as $n\to\infty,$ the above representation together with~\eqref{P1} (noting that $nV(\ep_n M_n)\le H_n (x-kM_n)(1+o(1))$ for $|x-kM_n|<\ep_n$) and~\eqref{P2} completes the proof of Theorem~\ref{thm_kM}.

\hfill$\Box$

\medskip

{\em Proof of Theorem~\ref{thm_kinM}.} We fix $k\ge 1$  and again start with~\eqref{1_sum}, but now with a new, much lower threshold~$y:=h_n M_n$, where $h_n\downarrow 0$ slowly enough so that   $h_nM_n\gg s_n$ (and hence still $nV(h_n M_n)\to 0$ as $n\to\infty$) and
\begin{align}
\label{hML}
L(h_n M_n) \sim  L(M_n).
\end{align}
That this is always possible can be demonstrated applying the standard argument to the claim of the uniform convergence theorem for slowly varying functions.

Using~\eqref{nu>k}, for the last term on the right-hand side of~\eqref{1_sum} one obtains
\begin{align}
\label{oPi}
\pr (S_n >x, \nu_n (h_n M_n) >k) =O\big((nV(h_nM_n))^{k+1}\big)=o\big(\Pi_n ^{k}\big)
\end{align}
when we take $h_n\ge n^{-(\al -2)/(4\al (k+1))}.$ Indeed, in this case, as $M_n\gg n^{1/2}$ and $L(M_n)=o(n^\delta)$ for any $\delta>0,$ one has from~\eqref{hML} that
\begin{align}
(nV(h_nM_n))^{k+1}
& \sim  n^{k+1} h_n^{-\al (k+1)} M_n^{-\al (k+1)} L(M_n)^{k+1}
\notag \\
& \le  (1+o(1)) \Pi_n ^k n^{\al/4 +1/2}   M_n^{-\al  } L(M_n)
\notag \\
& \ll \Pi_n ^k n^{-(\al -2)  /4}   L(M_n) =o(\Pi_n ^k).
\label{bound_hM}
\end{align}

To bound the terms in the first sum on the right-hand side of~\eqref{1_sum}, we note that, for $0\le j<k-1,$ one has $x-jM_n>M_n$. Then  we argue as in~\eqref{11_sum} to demonstrate, assuming without loss of generality that $h_n<1/(k+1),$ that 
\begin{align*}
 \pr (Y_n >x ,\nu_n (h_nM_n)=j)
 & \le  {n\choose j} V^j (h_nM_n) \pr \big(S_{n-j} >x-jM_n, \overline{\xi}_{n-j}\le h_nM_n   \big)
 \\
  & \le  {n\choose j} V^j (h_n M_n) \pr \big(S_{n-j} > M_n, \overline{\xi}_{n-j}\le M_n/(k+1)   \big)
  \\ & = O\big( (nV(h_n M_n))^j \Pi_n ^{k+1/2 }  \big)=o(\Pi_n ^k),
\end{align*}
where we used Corollary~4.1.3 from~\cite{BoBo08} in the second last relation.

For $j=k-1$ we obtain, for $x>(k-1+h)M_n,$ $h_n\le 1/3,$
\begin{align*}
 \pr (Y_n >x & ,  \nu_n (h_n M_n)  =k-1)\\
 &
 \le  {n\choose k-1} V^{k-1} (h_n M_n) \pr \big(S_{n-j} > h_n M_n; \overline{\xi}_{n-k+1} \le h_n M_n  \big)
 \\
 & =   O\big( (nV(h_nM_n))^{k-1} \big) \cdot \pr \big(S_{n-j} > h_n M_n; \overline{\xi}_{n-k+1}\le M_n/3  \big)
 \\
 & =  O\big( (nV(h_n M_n))^{k+1} \big) =o(\Pi_n ^k),
\end{align*}
again using Corollary~4.1.3 from~\cite{BoBo08} and~\eqref{bound_hM}.

Combining~\eqref{oPi} with the last two  bounds and Remark~\ref{rem_1}, we conclude that it remains to demonstrate that the asymptotics claimed in~\eqref{k_sum} hold true for the middle term on the right-hand side of~\eqref{1_sum} with $y=h_n M_n,$ i.e.\ for
\begin{align}
\pr (Y_n >x & , \nu_n (h_n M_n)  =k)
    =  {n \choose k} \pr \big(Y_n >x, \overline{\xi}_{n-k} \le h_n M_n <  \underline{\xi}_{n,k}\big)
  \notag \\
  & =  {n \choose k} \pr \big(S_{n-k} + Y_{n,k}>x, \overline{\xi}_{n-k} \le h_n M_n <  \underline{\xi}_{n,k}\big)
   \notag \\
  & =  {n \choose k} \big[\pr \big(B;   \underline{\xi}_{n,k}>h_n M_n \big)
   -\pr \big(B; \overline{\xi}_{n-k} \wedge   \underline{\xi}_{n,k}>h_n M_n \big)\big],
  \label{case_k}
\end{align}
where we put $Y_{n,k}:=\sum_{j=n-k+1}^n X_{n,j}$  and $B:=\{ S_{n-k}  + Y_{n,k}>x \}.$

First we will turn our attention to  the first term in the square brackets in~\eqref{case_k}.   Setting    $C:= \{|S_{n-k} |\le 2s_n\}$ and $x_\pm:= x\pm  2s_n,$ we see that 
\[
\pr \big(B;   \underline{\xi}_{n,k}>h_n M_n \big) = \pr \big(BC;   \underline{\xi}_{n,k}>h_n M_n \big) + \pr \big(BC^c;   \underline{\xi}_{n,k}>h_n M_n \big),
\]
where
\begin{align}
\pr \big(  Y_{n,k} & >x_+,   \underline{\xi}_{n,k}>  h_n M_n  \big)
    - \pr \big( C^c;  \underline{\xi}_{n,k}>  h_n M_n  \big)
\notag
\\
& \le   \pr \big(BC;   \underline{\xi}_{n,k}>  h_n M_n  \big)
   \le
   \pr \big(  Y_{n,k}>x_-,   \underline{\xi}_{n,k}>  h_nM_n  \big)
   \label{bounds_S}
\end{align}
and, using  Chebyshev's inequality and representing $nV(h_n M_n)$ as in the first line of~\eqref{bound_hM},
\begin{align*}
 {n\choose k} \pr \big(BC^c;   \underline{\xi}_{n,k}>  h_n M_n   \big)
 & \le {n\choose k} \pr \big( C^c;   \underline{\xi}_{n,k}>  h_n M_n   \big)
 \\
  & \le
  \frac1{k!} \pr \big(|S_{n-k}|> 2s_n) (nV (h_n M_n) )^k
  \\
  & =O\Big( \frac{n}{n \ln n} \cdot  h_n^{-\al k}\Pi_n ^k  \Big) =o (\Pi_n ^k)
\end{align*}
once we choose $h\gg (\ln n)^{-1/(\al k)}.$

To bound the second term in the square brackets in~\eqref{case_k},  we observe that
\begin{align*}
  {n\choose k}\pr \big(B; \overline{\xi}_{n-k} \wedge   \underline{\xi}_{n,k}>h_n M_n \big)
  & \le  \frac{n^k}{k!} \pr \big(  \overline{\xi}_{n-k} > h_n M_n \big)
  \pr \big(    \underline{\xi}_{n,k}>h_n M_n\big)
  \\
   & \le \frac{n}{k!}   V(h_n M_n) \cdot  (nV(h_n M_n))^k
   \\
   & =O\big((nV(h_n M_n))^{k+1}\big)
    =o (\Pi_n ^k)
\end{align*}
from~\eqref{bound_hM}.

The above bounds show that the expression in the last line of~\eqref{case_k}   is ``squeezed" between the values $ {n\choose k}\pr \big(  Y_{n,k}   >x_\pm,   \underline{\xi}_{n,k}>  h_n M_n  \big)+ o (\Pi_n ^k).$ Therefore, to complete the proof of the theorem, it suffices to  demonstrate that, for $z=x_\pm,$  the expression
\begin{align}
 {n \choose k}\pr \big(  Y_{n,k}>z ,   \underline{\xi}_{n,k}>  h_n M_n  \big)
 = {n \choose k}\pr \big(  Y_{ k}>z ,   \underline{\xi}_{ k}>  h_n M_n  \big)
\label{chooseS}
\end{align}
follows the asymptotics on the right-hand side of~\eqref{k_sum}.

To this end, first observe that,  since $\underline{\xi}_{k,j}$ is independent of $\xi_1,\ldots, \xi_{k-j},$ $0\le j \le k,$  for 
\[
z\in I_{k,n}:= ((k-1+h_n)M_n, (k-h_n)M_n) 
\]
one has
\begin{align}
\pr \big(  Y_{k} & >z ,   \underline{\xi}_{k}>  h_n M_n  \big)
  = \sum_{j=0}^k  \pr \big(  Y_{ k}>z ,   \underline{\xi}_{ k}>  h_n M_n, \nu_k (M_n) =j  \big)
  \notag\\
    &
  = \sum_{j=0}^k {k \choose j}\pr \big(   S_{ k-j}>z-j M_n ,
   h_n M_n < \underline{\xi}_{ k-j}, \overline{\xi}_{ k-j} \le M_n <\underline{\xi}_{k,j}  \big)
  \notag\\
    &
  = \sum_{j=0}^k {k \choose j}\pr \big(   S_{ k-j}>z-jM_n ,
  \overline{\xi}_{k-j} \le M_n , \underline{\xi}_{k-j} > h_n   M_n  \big)V^j (M_n).
  \label{Ps>z}
\end{align}
Clearly, the term in the last sum with $j=k$ is equal to just $V^k(M_n).$ To derive the asymptotic behavior of all the other terms in that sum,  note that, for any fixed $m\ge 1 $ and $\boldsymbol{u}:=(u_1, \ldots, u_m)\in [1,\infty)^m,$ due to~\eqref{regular_V} one has
\begin{align}
\pr \Bigl(\frac{\xi_j}{N}> u_j, j\in [m]
\Big| \underline{\xi}_m > N\Bigr)
 =
 \prod_{j\in [m]}\frac{V(u_j N)}{V( N)}
 \to \prod_{j\in [m]}u_j^{-\al}, \quad N\to \infty.
  \label{conv_cond}
\end{align}
Hence the conditional distribution of the random vector $(\xi_1, \ldots, \xi_m) /N$ given $\underline{\xi}_m > N$ converges weakly as $N\to\infty$ to the probability   law $\Theta_m$ with density $\al^m \prod_{j\in [m]}u_j^{-\al -1}$ on $[1,\infty)^m.$

Therefore, setting
\[
D_m(v,g):=\{\boldsymbol{u}\in \R^m :1<u_i<g, i\in [m], u_1+\cdots +u_m>v\},\quad v,g\ge 1,
\]
one has, due to the absolute continuity of the limiting distribution $\Theta_m$ and positivity of its density on $[1,\infty)^m$, the uniform relative convergence of the conditional distribution tail  of the scaled sum $S_m$: for any fixed  $h_0\in (0,m/2),$ as $n\to\infty,$
\begin{align}
\sup_{1\le v\le m/h_0-1}
\left|\frac{\pr \Big(   \frac{S_{m}}{h_0 M_n}>v , \frac{\overline{\xi}_m}{h_0 M_n}\le \frac1{h_0} \, \Big|\, \underline{\xi}_m > h_0 M_n \Big)}{\Theta_m(D_m(v,1/h_0))}-1\right| \to 0.
\label{unif_con}
\end{align}
The standard argument shows that \eqref{unif_con} will still hold if one replaces in it the fixed $h_0>0$ with our positive sequence $h_n\downarrow 0$ provided that the latter is vanishing  slowly enough as $n\to \infty.$  

Now observe that, after the above-mentioned replacement of $h_0$ with $h_n$, the probability of the thus  modified condition $\{\underline{\xi}_m > h_n M_n\}$ from~\eqref{unif_con} is equal, in view of~\eqref{hML},  to
\begin{align}
V^m(h_n  M_n) & =(h_n M_n)^{-\al m} L^m (h_n M_n)
 \notag \\
& \sim  (h_n M_n)^{-\al m} L^m (M_n)
\sim   h_n  ^{-\al m} V^m (M_n) .
\label{VhM}
\end{align}
We conclude that
\begin{align}
\sup_{1\le v\le m/h_n  -1}
\bigg|\frac{\pr \big(    S_{m} >v h_n   M_n ,  \overline{\xi}_m \le M_n , \underline{\xi}_m > h_n   M_n \big)}{h_n  ^{-\al m} V^m (M_n) \Theta_m(D_m(v,1/h_n  ))}-1\bigg| \to 0.
\label{unif}
\end{align}
Changing the variables $u_i:=t_i/h_n ,$  $i\in [m],$ we see that
\begin{align}
 h_n  ^{-\al m} \Theta_m(D_m(v &,1/h_n  ))  = \al^m \int_{D_m(v,1/h_n )}h_n  ^{-\al m}  (u_1\cdots u_m)^{-\al -1}d\boldsymbol{u}
 \notag \\
 &
  =  \al^m \int_{D_m(vh_n ),\ \underline{t}_m>h_n  }  (t_1\cdots t_m)^{-\al -1}d\boldsymbol{t}
  = W_m (vh_n )
  \label{Wm}
\end{align}
provided that $vh_n  >m-1 +h_n ,$ since  then automatically $\underline{t}_m>h_n $ for any $\boldsymbol{t}\in D_m(vh_n ).$

Now clearly $(z-jM_n)/M_n> k-j -1 +h_n \ge h_n $ for $z\in I_{k,n},$ $0\le j<k,$ so that, due to~\eqref{unif} and~\eqref{Wm} (where we put $vh_n M_n=z-jM_n$, $m=k-j$), one has
\begin{align}
\sup_{z\in I_{k,n}}
\bigg| \frac{\pr \big(   S_{ k-j}>z-jM_n ,
  \overline{\xi}_{k-j} \le M_n , \underline{\xi}_{k-j} > h_n   M_n  \big)}{W_{k-j} (z/M_n-j)V^{k-j}( M_n)}
  -1\bigg| \to 0,
  \label{ratioW}
\end{align}
the convergence being uniform in $z\in I_{k,n}.$

Note that  ${n \choose k}=(1+\theta_{n,k})n^k/k!,$ where $\theta_{n,k}:=\prod_{i=1}^{k-1}(1-i/n)-1\to 0 $ as $n\to\infty.$ Hence it follows from~\eqref{Ps>z} and~\eqref{ratioW} that the expression on the right-hand side of~\eqref{chooseS} is equal to
\begin{align*}
\frac{(1+\theta_{n,k})\Pi_n ^k  }{k!}
&
\sum_{j=0}^k {k \choose j}\pr \big(   S_{ k-j}>z-jM_n ,  \overline{\xi}_{k-j} \le M_n , \underline{\xi}_{k-j} > h_n M_n     \big) V^{j-k}(M_n)
\\
 & = \frac{(1+o(1))\Pi_n ^k  }{k!}
\sum_{j=0}^k {k \choose j}W_{k-j} (z/M_n-j) ,
\end{align*}
the term $o(1)$ being uniform in $z\in I_{k,n}.$   Since clearly $x_\pm = x\pm 2s_n \in I_{k,n}$ for $x\in ((k -1  + h_{k,n})M_n,(k - h_{k,n})M_n) $ with $h_{k,n}:= h_n+2s_n/M_n\sim h_n ,$ the desired behavior of~\eqref{chooseS} is established. Theorem~\ref{thm_kinM} is proved.

\hfill $\Box$

\end{document}